\documentclass[11pt,reqno]{article}
\usepackage{amsmath,mathrsfs,graphicx,amssymb,amsfonts,color}
\usepackage{cite}

\font\bg=cmbx10 scaled\magstep1

\font\small=cmr8

\newtheorem{newlemma}{{\bf Lemma}}

\newenvironment{lema}{\begin{newlemma}{\hspace{-0.5
em}{\bf.}}}{\end{newlemma}}

\newtheorem{newteorem}{{\bf Theorem}}

\newenvironment{teorem}{\begin{newteorem}{\hspace{-0.5
em}{\bf.}}}{\end{newteorem}}

\newtheorem{newkorolari}{{\bf Corollary}}

\newenvironment{korolari}{\begin{newkorolari}{\hspace{-0.5
em}{\bf.}}}{\end{newkorolari}}

\newtheorem{newdefine}{{\bf Definition}}

\newtheorem{newquestion}{{\bf Question}}

\newtheorem{newkonjek}{{\bf Conjecture}}

\newtheorem{newexample}{{\bf Example}}

\begin{document}
\tolerance=10000
\baselineskip18truept
\newbox\thebox
\global\setbox\thebox=\vbox to 0.2truecm{\hsize 0.15truecm\noindent\hfill}
\def\boxit#1{\vbox{\hrule\hbox{\vrule\kern0pt
\vbox{\kern0pt#1\kern0pt}\kern0pt\vrule}\hrule}}
\def\qed{\lower0.1cm\hbox{\noindent \boxit{\copy\thebox}}\bigskip}
\def\nt{\noindent}

\centerline{\Large \bf Randi\'{c} energy of specific  graphs }
\vspace{.3cm}

\bigskip

\baselineskip12truept
\centerline{Saeid  Alikhani$^{}${}\footnote{\baselineskip12truept\it\small
Corresponding author. E-mail: alikhani@yazd.ac.ir}and  Nima Ghanbari }
\baselineskip20truept
\centerline{\it Department of Mathematics, Yazd University}
\vskip-8truept
\centerline{\it  89195-741, Yazd, Iran}

\vskip-0.2truecm
\noindent\rule{16cm}{0.1mm}
\noindent{\bg{Abstract}}

\baselineskip14truept

{\nt Let $G$ be a simple graph with vertex set $V(G) = \{v_1, v_2,\ldots, v_n\}$. The Randi\'{c} matrix of $G$, denoted by $R(G)$, is defined as the $n\times n$ matrix whose $(i,j)$-entry is $(d_id_j)^{\frac{-1}{2}}$ if $v_i$ and
$v_j$ are adjacent and $0$ for another cases.
Let the eigenvalues of the Randi\'{c} matrix $R(G)$ be  $\rho_1\geq \rho_2\geq \ldots\geq \rho_n$ which are the roots of the Randi\'c characteristic polynomial $\prod_{i=1}^n (\rho-\rho_i)$. The Randi\'{c} energy $RE$ of $G$ is the sum of
absolute values of the eigenvalues of $R(G)$. In this paper we compute the Randi\'c characteristic polynomial and the Randi\'c energy for  specific graphs $G$. 
}


\noindent\rule{16cm}{0.1mm}

\baselineskip20truept

\section{Introduction}

\nt In this paper we are concerned with simple finite graphs, without directed, multiple, or weighted edges, and without self-loops. Let $G$ be such a graph, with vertex set $V(G) = \{v_1, v_2,\ldots, v_n\}$. If two vertices $v_i$ and $v_j$ of $G$ are adjacent, then we use the notation $v_i \sim v_j $. For $v_i \in V(G)$, the degree of the vertex $v_i$, denoted by $d_i$, is the number of the vertices adjacent to $v_i$.

\nt  Let $A(G)$ be adjacency matrix of $G$ and $\lambda_1,\lambda_2,\ldots,\lambda_n$ its eigenvalues. These are said to be be the eigenvalues of the graph $G$ and to form its spectrum \cite{Cve}. The energy $E(G)$ of the graph $G$ is defined as the sum of the absolute values of its eigenvalues
$$E(G)=\sum_{i=1}^n\vert\lambda_i\vert.$$
Details and more information on graph energy can be found in \cite{Gut,Gut1,Gut2,Maj}.

\nt In 1975 Milan Randi\'{c} invented a molecular structure descriptor defined as \cite{Ran}
$$R(G)=\sum_{v_i\sim v_j} \frac{1}{\sqrt{d_i d_j}}.$$

\nt The Randi\'{c}-index-concept suggests that it is purposeful to associate to the graph $G$ a symmetric square matrix $R(G)$. The Randi\'{c} matrix $R(G)=(r_{ij})_{n\times n}$ is defined as \cite{Boz,Boz1,Gut3}
\begin{displaymath}
 r_{ij}= \left\{ \begin{array}{ll}
\frac{1}{\sqrt{d_i d_j}} & \textrm{if $v_i \sim v_j$}\\
0 & \textrm{otherwise.}
\end{array} \right.
\end{displaymath}

\nt Denote the eigenvalues of the Randi\'{c} matrix $R(G)$ by $\rho_1,\rho_2,\ldots,\rho_n$ and label them in non-increasing order. Similar to characteristic polynomial of a matrix, we consider  the Randi\'c characteristic polynomial of $R(G)$ (or a  graph $G$), as $det(\rho I-R(G))$ which is equal to $\prod_{i=1}^n (\rho-\rho_i)$. The Randi\'{c} energy \cite{Boz,Boz1,Gut3} of $G$ is defined as
$$RE(G)=\sum_{i=1}^n\vert\rho_i\vert.$$
For several lower and upper bounds on Randi\'{c} energy, see \cite{Boz,Boz1,Gut3}.

\nt  In Section 2, we obtain the Randi\'{c} characteristic polynomial and energy of specific graphs. As a result, we show that for every natural number $m\geq 2$, there exists a graph $G$ such that $RE(G)=m$. In Section 3, we find Randi\'{c} energy of specific graphs with one edge deleted.

\nt As usual we denote by $J$ a matrix all  of whose entries are equal to $1$.

\section{Randi\'{c} characteristic polynomial and Randi\'c energy of specific graphs }

\nt In this section we study the Randi\'c characteristic polynomial and the Randi\'{c} energy for certain graphs. The following theorem gives a relationship between the Randi\'c energy and energy of path $P_n$.


\begin{newlemma}\rm\cite{Gut3} Let $P_n$ be the path on $n$ vertices. Then
$$RE(P_n)=2+\frac{1}{2}E(P_{n-2}).$$
\end{newlemma}

\nt The following theorem gives the Randi\'c energy of even cycles.

\begin{newlemma}\rm\cite{Roj} Let $C_{2n}$ be the cycle on $2n$ vertices for $n\geq 2$. Then
$$RE(C_{2n})=\frac{2sin((\lfloor\frac{n}{2}\rfloor+\frac{1}{2})\frac{\pi}{n})}{sin\frac{\pi}{2n}}.$$
\end{newlemma}

\nt Here we shall compute the Randi\'c characteristic polynomial of paths and cycles.

\begin{teorem}\label{path}
For $n\geq 5$, the Randi\'{c} characteristic polynomial of the path graph $P_n$ satisfy
$$RP(P_n,\lambda)=(\lambda^2-1)(\lambda \Lambda_{n-3}-\frac{1}{4}\Lambda_{n-4}),$$
where for every $k\geq 3$, $\Lambda_k=\lambda \Lambda_{k-1}-\frac{1}{4}\Lambda_{k-2}$ with  $\Lambda_1=\lambda$ and $\Lambda_2=\lambda ^2-\frac{1}{4}$.
\end{teorem}

\noindent{\bf Proof.}
\nt For every $k\geq 3$, consider

$$B_k :=
\left(\begin{array}{ccccccccc}
\lambda& \frac{-1}{2}&0&0&\ldots &0&0 &0 \\
 \frac{-1}{2}& \lambda  &\frac{-1}{2} &0&\ldots &0&0  &0 \\
0&  \frac{-1}{2}& \lambda &\frac{-1}{2} & \ldots &0&0  &0 \\
0&0&  \frac{-1}{2}& \lambda  & \ldots &0&0  &0 \\
\vdots & \vdots & \vdots & \vdots &\ddots &\vdots&\vdots &\vdots \\
 0& 0&0&0  &\ldots &\lambda&\frac{-1}{2} &0 \\
0&  0&0&0  &\ldots &\frac{-1}{2}&\lambda&\frac{-1}{2}  \\
0&  0& 0&0 &\ldots &0&\frac{-1}{2}&\lambda   \\
\end{array}\right)_{k\times k}, $$

\nt and let $\Lambda_k=det(B_k)$. It is easy to see that $\Lambda_k=\lambda \Lambda_{k-1}-\frac{1}{4}\Lambda_{k-2}$.

\nt Suppose that  $RP(P_n,\lambda)=det(\lambda I - R(P_n) )$. We have

$$ RP(P_n,\lambda) = det
\left(\begin{array}{c|ccccc|c}
\lambda& \frac{-1}{\sqrt{2}}&0&\ldots &0&0 &0 \\
\hline
 \frac{-1}{\sqrt{2}}& &&&& &0 \\
0&  &&&&  &0 \\
\vdots & &&B_{n-1}&&&\vdots \\
 0&&&&& &0 \\
0&  &&&&&\frac{-1}{\sqrt{2}}  \\
\hline
0&  0&0 &\ldots &0&\frac{-1}{\sqrt{2}}&\lambda   \\
\end{array}\right)_{n\times n}.$$

\nt So

$$ RP(P_n,\lambda) =\lambda
det\left(\begin{array}{cccc|c}
&&&  &0 \\
&&&&\vdots \\
 &&B_{n-2}& &0 \\
 &&&&\frac{-1}{\sqrt{2}}  \\
\hline
0  &\ldots &0&\frac{-1}{\sqrt{2}}&\lambda   \\
\end{array}\right)+
\frac{1}{\sqrt{2}}det\left(\begin{array}{c|ccc|c}
\frac{-1}{\sqrt{2}}& \frac{-1}{2}&\ldots &0 &0 \\
\hline
0& && &0 \\
\vdots & &B_{n-3}&&\vdots \\
0&  &&&\frac{-1}{\sqrt{2}}  \\
\hline
0&  0 &\ldots &\frac{-1}{\sqrt{2}}&\lambda   \\
\end{array}\right).
$$

\nt And so

$$ RP(P_n,\lambda) =\lambda (\lambda \Lambda_{n-2}+
\frac{1}{\sqrt{2}}det\left(\begin{array}{ccc|c}
 &&  &0\\
&B_{n-3}&&\vdots  \\
 &&&\frac{-1}{2}  \\
\hline
0 &\ldots &0&\frac{-1}{\sqrt{2}}   \\
\end{array}\right))-
\frac{1}{2}det\left(\begin{array}{ccc|c}
 &&  &0\\
&B_{n-3}&&\vdots  \\
 &&&\frac{-1}{\sqrt{2}} \\
\hline
0 &\ldots &\frac{-1}{\sqrt{2}}&\lambda   \\
\end{array}\right).
$$

\nt Therefore

$$ RP(P_n,\lambda) =\lambda^2  \Lambda_{n-2}-
\frac{1}{2}\lambda \Lambda_{n-3}-
\frac{1}{2}\lambda \Lambda_{n-3}-\frac{\sqrt{2}}{4}det\left(\begin{array}{ccc|c}
 &&  &0\\
&B_{n-4}&&\vdots  \\
 &&&\frac{-1}{2}  \\
\hline
0 &\ldots &0&\frac{-1}{\sqrt{2}}   \\
\end{array}\right),
$$

\nt and so
$$ RP(P_n,\lambda)=\lambda^2  \Lambda_{n-2}-\lambda \Lambda_{n-3}+\frac{1}{4} \Lambda_{n-4}=\lambda^2(\lambda \Lambda_{n-3}
-\frac{1}{4}\Lambda_{n-4})-\lambda \Lambda_{n-3}+\frac{1}{4}\Lambda_{n-4}.$$

\nt Hence

$$ RP(P_n,\lambda)=(\lambda^2-1)(\lambda \Lambda_{n-3}-\frac{1}{4}\Lambda_{n-4}).\quad\qed$$

\begin{teorem}
For $n\geq 3$, the Randi\'{c} characteristic polynomial of the cycle graph $C_n$ is
$$RP(C_n,\lambda)=\lambda \Lambda_{n-1}-\frac{1}{2}\Lambda_{n-2}-(\frac{1}{2})^{n-1},$$
where for every $k\geq 3$, $\Lambda_k=\lambda \Lambda_{k-1}-\frac{1}{4}\Lambda_{k-2}$ with  $\Lambda_1=\lambda$ and $\Lambda_2=\lambda ^2-\frac{1}{4}$.
\end{teorem}
\noindent{\bf Proof.}  Similar to the proof of Theorem \ref{path}, for every $k\geq 3$, we consider

$$B_k :=
\left(\begin{array}{ccccccccc}
\lambda& \frac{-1}{2}&0&0&\ldots &0&0 &0 \\
 \frac{-1}{2}& \lambda  &\frac{-1}{2} &0&\ldots &0&0  &0 \\
0&  \frac{-1}{2}& \lambda &\frac{-1}{2} & \ldots &0&0  &0 \\
0&0&  \frac{-1}{2}& \lambda  & \ldots &0&0  &0 \\
\vdots & \vdots & \vdots & \vdots &\ddots &\vdots&\vdots &\vdots \\
 0& 0&0&0  &\ldots &\lambda&\frac{-1}{2} &0 \\
0&  0&0&0  &\ldots &\frac{-1}{2}&\lambda&\frac{-1}{2}  \\
0&  0& 0&0 &\ldots &0&\frac{-1}{2}&\lambda   \\
\end{array}\right)_{k\times k}, $$

\nt and let $\Lambda_k=det(B_k)$. We have $\Lambda_k=\lambda \Lambda_{k-1}-\frac{1}{4}\Lambda_{k-2}$.

\nt Suppose that $RP(C_n,\lambda)=det(\lambda I - R(C_n) )$. We have
$$  RP(C_n,\lambda) =det
\left(\begin{array}{c|ccccc}
\lambda& \frac{-1}{2}&0&\ldots &0&\frac{-1}{2} \\
\hline
 \frac{-1}{2}& &&&& \\
0&  &&&&  \\
\vdots & &&B_{n-1}&& \\
 0&&&&& \\
\frac{-1}{2}&  &&&&\\
\end{array}\right)_{n\times n}.$$

\nt So

$$  RP(C_n,\lambda) =\lambda \Lambda_{n-1}+\frac{1}{2}det
\left(\begin{array}{c|ccc}
\frac{-1}{2}& \frac{-1}{2}&\ldots &0 \\
\hline
0&  &&  \\
\vdots  &&B_{n-2}& \\
\frac{-1}{2}&  &&\\
\end{array}\right)+(-1)^{n+1}(\frac{-1}{2})det
\left(\begin{array}{c|ccc}
\frac{-1}{2}& && \\
\vdots  &&B_{n-2}& \\
0&  &&  \\
\hline
\frac{-1}{2}&0  &\ldots &\frac{-1}{2}\\
\end{array}\right).
$$

\nt And so,

$  RP(C_n,\lambda) =\lambda \Lambda_{n-1}-\frac{1}{4}\Lambda_{n-2}+(-1)^n(\frac{-1}{4})det
\left(\begin{array}{ccc|c}
\frac{-1}{2} &\ldots&0  &0\\
\hline
&&&0 \\
& B_{n-3}&& \vdots \\
&  &&\frac{-1}{2}\\
\end{array}\right)+\\
(-1)^{n+1}(\frac{-1}{2})(\frac{-1}{2}det
\left(\begin{array}{c|ccc}
\frac{-1}{2}& && \\
\vdots  &&B_{n-3}& \\
0&  &&  \\
\hline
\frac{-1}{2}&0  &\ldots &\frac{-1}{2}\\
\end{array}\right)+(-1)^n(\frac{-1}{2})\Lambda_{n-2}
).
$

\nt Therefore,
$$RP(C_n,\lambda) =\lambda \Lambda_{n-1}-\frac{1}{4}\Lambda_{n-2}+(-1)^n(\frac{-1}{2})^{n-1}(\frac{1}{2})
+(-1)^{n+1}(\frac{-1}{2})^n+\frac{1}{4}(-1)^{2n+1}\Lambda_{n-2}.
$$

\nt Hence

$$RP(C_n,\lambda) =\lambda \Lambda_{n-1}-\frac{1}{2}\Lambda_{n-2}-(\frac{1}{2})^{n-1}.\quad\qed
$$

\nt In \cite{Roj} has shown that the Randi\'{c} energy of $S_n=K_{1,n-1}$, the star on $n$ vertices and the complete bipartite graph $K_{m,n}$ is $2$. Here using the Randi\'c characteristic polynomial, we prove these results. We need the following lemma: 

\begin{lema} \label{new}\rm\cite{Cve}
If $M$ is a nonsingular square matrix, then
$$det\left(  \begin{array}{cc}
M&N  \\
P& Q \\
\end{array}\right)=det (M) det( Q-PM^{-1}N).
$$
\end{lema}

\begin{teorem}
For $n\geq 2$,
\begin{itemize}
\item[(i)] The Randi\'{c} characteristic polynomial of the star graph $S_n=K_{1,n-1}$ is
$$RP(S_n,\lambda)=\lambda^{n-2}(\lambda ^2 -1).$$
\item[(ii)] The Randi\'{c} energy of $S_n$ is
$$RE(S_n)=2.$$
\end{itemize}
\end{teorem}
\noindent{\bf Proof.}
\begin{enumerate}
\item[(i)] It is easy to see that the Randi\'c matrix of $K_{1,n-1}$ is
$\frac{1}{\sqrt{n-1}}\left( \begin{array}{cc}
0_{1\times 1}&J_{1\times {n-1}} \\
J_{{n-1}\times 1}&0_{{n-1}\times {n-1}}   \\
\end{array} \right)$.
We  have 
$$det(\lambda I -R(S_n))=det
\left(  \begin{array}{cc}
\lambda  & \frac{-1}{\sqrt{n-1}}J_{1\times (n-1)}  \\
\frac{-1}{\sqrt{n-1}}J_{(n-1)\times 1}& \lambda I_{n-1} \\
\end{array}\right).
$$
Using Lemma \ref{new},
$$
det(\lambda I -R(S_n))=\lambda det( \lambda I_{n-1} - \frac{1}{\sqrt{n-1}}J_{(n-1)\times 1}\frac{1}{\lambda} \frac{1}{\sqrt{n-1}}J_{1\times (n-1)}).
$$

We know that $J_{(n-1)\times 1}J_{1\times (n-1)}=J_{n-1}$. Therefore
$$
det(\lambda I -R(S_n))=\lambda det( \lambda I_{n-1} - \frac{1}{\lambda (n-1)}J_{n-1})=\lambda ^{2-n}
det( \lambda ^2 I_{n-1} - \frac{1}{n-1}J_{n-1}).
$$
Since the eigenvalues of $J_{n-1}$ are $n-1$ (once) and 0 ($n-2$ times),  the eigenvalues of $\frac{1}{n-1}J_n$ are $1$ (once) and 0 ($n-2$ times). Hence
 $$RP(S_n,\lambda)=\lambda^{n-2}(\lambda ^2 -1).$$

 \item[(ii)] It follows from Part (i).\quad\qed

 \end{enumerate}

\begin{teorem}
For $n\geq 2$,
\begin{itemize}
\item[(i)] the Randi\'{c} characteristic polynomial of complete graph $K_n$ is
$$RP(K_n,\lambda)=(\lambda -1)(\lambda +\frac{1}{n-1})^{n-1}.$$
\item[(ii)] the Randi\'{c} energy of $K_n$ is
$$RE(K_n)=2.$$
\end{itemize}
\end{teorem}
\noindent{\bf Proof.}
\begin{enumerate}
\item[(i)]
It is easy to see that the Randi\'{c} matrix of $K_n$ is $\frac{1}{n-1}(J-I)$. Therefore
$$
RP(K_n,\lambda)=det(\lambda I- \frac{1}{n-1}J+ \frac{1}{n-1}I)=det((\lambda +\frac{1}{n-1})I- \frac{1}{n-1}J).
$$

Since the eigenvalues of $J_n$ are $n$ (once) and 0 ($n-1$ times),  the eigenvalues of $\frac{1}{n-1}J_n$ are $\frac{n}{n-1}$ (once) and 0 ($n-1$ times). Hence
$$RP(K_n,\lambda)=(\lambda -1)(\lambda +\frac{1}{n-1})^{n-1}.$$

\item[(ii)] It follows from Part (i).\quad\qed
\end{enumerate}

\begin{teorem}\label{bipartite}
For natural number $m,n\neq 1$,
\begin{itemize}
\item[(i)] The Randi\'{c} characteristic polynomial of complete bipartite graph $K_{m,n}$ is
$$RP(K_{m,n},\lambda)=\lambda^{m+n-2}(\lambda^2 -1).$$
\item[(ii)] The Randi\'{c} energy of $K_{m,n}$ is
$$RE(K_{m,n})=2.$$
\end{itemize}
\end{teorem}
\noindent{\bf Proof.}
\begin{enumerate}
\item[(i)]
It ie easy to see that the Randi\'{c} matrix of $K_{m,n}$ is
$\frac{1}{\sqrt{mn}}\left( \begin{array}{cc}
0_{m\times m}&J_{m\times n} \\
J_{n\times m}&0_{n\times n}   \\
\end{array} \right)$.
Using Lemma \ref{new} we have 

$$det(\lambda I -R(K_{m,n}))=det
\left(  \begin{array}{cc}
\lambda I_m & \frac{-1}{\sqrt{mn}}J_{m\times n}  \\
\frac{-1}{\sqrt{mn}}J_{n\times m}& \lambda I_n \\
\end{array}\right).
$$
So
$$
det(\lambda I -R(K_{m,n}))=det (\lambda I_m) det( \lambda I_n - \frac{1}{\sqrt{mn}}J_{n\times m}\frac{1}{\lambda}I_m \frac{1}{\sqrt{mn}}J_{m\times n}).
$$

We know that $J_{n\times m}J_{m\times n}=mJ_n$. Therefore
$$
det(\lambda I -R(K_{m,n}))=\lambda ^m det( \lambda I_n - \frac{1}{\lambda n}J_n)=\lambda ^{m-n}
det( \lambda ^2 I_n - \frac{1}{n}J_n).
$$
The eigenvalues of $J_n$ are $n$ (once) and 0 ($n-1$ times). So the eigenvalues of $\frac{1}{n}J_n$ are $1$ (once) and 0 ($n-1$ times).
Hence
$$RP(K_{m,n},\lambda)=\lambda^{m+n-2}(\lambda^2 -1).$$

\item[(ii)] It follows from Part (i).\quad\qed
\end{enumerate}

\medskip

\begin{figure}[!h]
\hspace{3.6cm}
\includegraphics[width=8.5cm,height=2.cm]{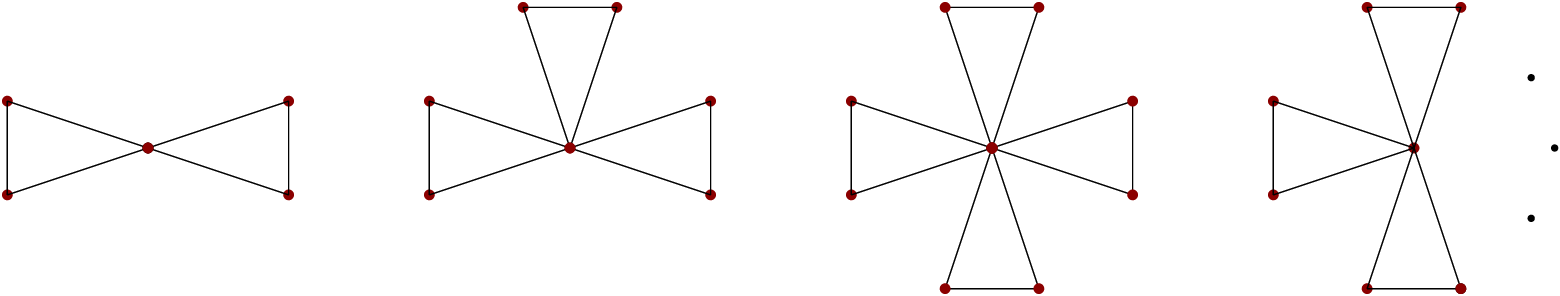}
\caption{\label{friend} Friendship graphs $F_2, F_3, F_4$ and $F_n$, respectively. }
\end{figure}

\nt Let $n$ be any positive integer and  $F_n$ be friendship graph with
$2n + 1$ vertices and $3n$ edges. In other words, the friendship  graph $F_n$ is a graph that can be constructed by coalescence $n$
copies of the cycle graph $C_3$ of length $3$ with a common vertex. The Friendship Theorem of  Erd\H os,
R\'{e}nyi and S\'{o}s \cite{erdos}, states that graphs with the property that every two vertices have
exactly one neighbour in common are exactly the friendship graphs.
The Figure \ref{friend} shows some examples of friendship graphs. Here we shall investigate the Randi\'{c} energy of friendship graphs.

\begin{teorem}
For $n\geq 2$,
\begin{itemize}
\item[(i)] The Randi\'{c} characteristic polynomial of friendship graph $F_n$ is
$$RP(F_n,\lambda)=(\lambda^2-\frac{1}{4})^{n-1}(\lambda -1)(\lambda +\frac{1}{2})^2.$$
\item[(ii)] The Randi\'{c} energy of $F_n$ is
$$RE(F_n)=n+1.$$
\end{itemize}
\end{teorem}
\noindent{\bf Proof.}
\begin{enumerate}
\item[(i)]
The Randi\'{c} matrix of $F_n$ is

$$
R(F_n) =
\left( \begin{array}{cccccc}
0 & \frac{1}{2\sqrt{n}} &\frac{1}{2\sqrt{n}}&\cdots & \frac{1}{2\sqrt{n}}& \frac{1}{2\sqrt{n}}\\
\frac{1}{2\sqrt{n}}& 0& \frac{1}{2}&\ldots &0&0  \\
\frac{1}{2\sqrt{n}}& \frac{1}{2}& 0&\ldots &0&0   \\
\vdots & \vdots &\vdots &\ddots &\vdots  \\
\frac{1}{2\sqrt{n}} & 0&0 &\ldots &0&\frac{1}{2}  \\
\frac{1}{2\sqrt{n}} & 0& 0&\ldots &\frac{1}{2}&0   \\
\end{array} \right)_{(2n+1)\times (2n+1)}.
$$

\nt Now for computing $det(\lambda I - R(F_n) )$, we consider its first row. The cofactor of the first array in this row is

\[\left(\begin{array}{cccccc}
\lambda& \frac{-1}{2}&\ldots &0&0  \\
  \frac{-1}{2}& \lambda&\ldots &0&0   \\
\vdots & \vdots &\ddots &\vdots&\vdots  \\
  0&0 &\ldots &\lambda&\frac{-1}{2}  \\
  0& 0&\ldots &\frac{-1}{2}&\lambda   \\
\end{array}\right) \]
\nt and the cofactor of  another arrays in the first row are similar to

\[
\left( \begin{array}{cccccc}
\frac{-1}{2\sqrt{n}}& \frac{-1}{2}&\ldots &0&0  \\
\frac{-1}{2\sqrt{n}}& \lambda&\ldots &0&0   \\
\vdots & \vdots &\ddots &\vdots&\vdots  \\
\frac{-1}{2\sqrt{n}}&0 &\ldots &\lambda&\frac{-1}{2}  \\
 \frac{-1}{2\sqrt{n}}& 0&\ldots &\frac{-1}{2}&\lambda   \\
\end{array} \right)
\]

\nt Now, by straightforward computation we have the result.
	
\item[(ii)] It follows from Part (i). \quad\qed

\end{enumerate}

\noindent{\bf Remark.} In \cite{Bap} has shown that the energy of a graph cannot be an odd integer. Since $RE(F_n)=n+1$ for $n\geq 2$, the Randi\'{c} energy  can be odd or even integer. More precisely we have:

\begin{korolari}
For every natural number $m\geq 2$, there exists a graph $G$ such that $RE(G)=m$.
\end{korolari}
\noindent{\bf Proof.} If $m=2$ then consider $K_2$ and for $m\geq 3$  we consider friendship graphs  $F_{m-1}$.\quad\qed

\begin{figure}[!h]
\hspace{1.6cm}
\includegraphics[width=11.5cm,height=4cm]{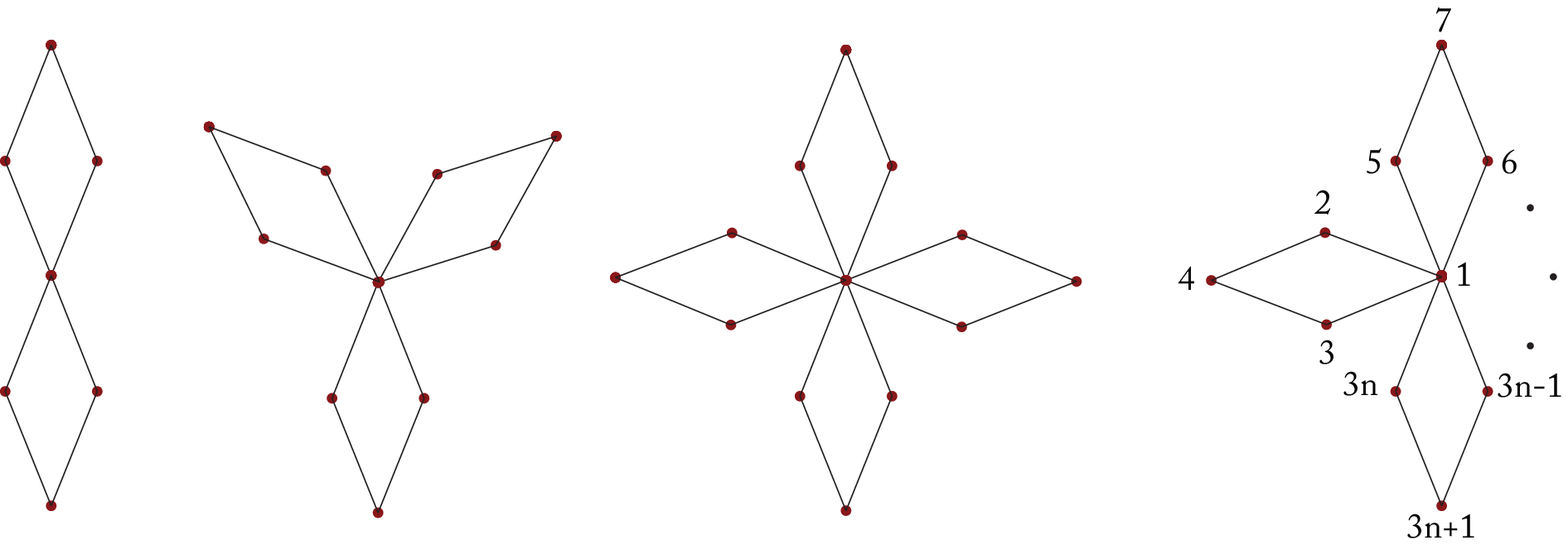}
\caption{\label{Dutch} Dutch Windmill Graph $D_4 ^2, D_4^3, D_4^4$ and $D_4^n$, respectively. }
\end{figure}

\nt Let $n$ be any positive integer and  $D_4^n$ be Dutch Windmill Graph with
$3n + 1$ vertices and $4n$ edges. In other words, the  graph $D_4^n$ is a graph that can be constructed by coalescence $n$
copies of the cycle graph $C_4$ of length $4$ with a common vertex.
The Figure \ref{Dutch} shows some examples of Dutch Windmill graphs. Here we shall investigate the Randi\'{c} energy of Dutch Windmill graphs.

\begin{teorem}
For $n\geq 2$,
\begin{itemize}
\item[(i)] The Randi\'{c} characteristic polynomial of friendship graph $D_4^n$ is
$$RP(D_4^n,\lambda)=\lambda^{n+1}(\lambda^2 -\frac{1}{2})^{n-1}(\lambda^2 -1).$$
\item[(ii)] The Randi\'{c} energy of $F_n$ is
$$RE(D_4^n)=2+(n-1)\sqrt{2}.$$
\end{itemize}
\end{teorem}
\noindent{\bf Proof.}
\begin{enumerate}
\item[(i)]
The Randi\'{c} matrix of $D_4^n$ is

$$
R(D_4^n) =
\left( \begin{array}{cccccccc}
0 & \frac{1}{2\sqrt{n}} &\frac{1}{2\sqrt{n}}&0&\cdots & \frac{1}{2\sqrt{n}}& \frac{1}{2\sqrt{n}}&0\\
\frac{1}{2\sqrt{n}}&0& 0& \frac{1}{2}&\ldots &0 &0&0 \\
\frac{1}{2\sqrt{n}}&0& 0& \frac{1}{2}&\ldots &0  &0&0\\
0&\frac{1}{2}& \frac{1}{2}&0 &\ldots &0&0 &0\\
\vdots & \vdots& \vdots &\vdots &\ddots &\vdots &\vdots \\
\frac{1}{2\sqrt{n}} & 0&0 &0&\ldots &0&0&\frac{1}{2}  \\
\frac{1}{2\sqrt{n}} & 0&0 &0&\ldots &0 &0&\frac{1}{2}  \\
0 & 0& 0&0&\ldots &\frac{1}{2}&\frac{1}{2}&0 \\
\end{array} \right)_{(3n+1)\times (3n+1)}.
$$

\nt Let $A= \left( \begin{array}{ccc}
 \lambda &0 &\frac{-1}{2}  \\
 0&\lambda &\frac{-1}{2} \\
\frac{-1}{2} &\frac{-1}{2} & \lambda \\
\end{array}\right)$,
$B= \left( \begin{array}{ccc}
 \frac{-1}{2\sqrt{n}} &0 &0  \\
\frac{-1}{2\sqrt{n}}&0 &0 \\
0 &0 & 0 \\
\end{array}\right)$
and
$C= \left( \begin{array}{ccc}
 \frac{-1}{2\sqrt{n}} & 0 &\frac{-1}{2}  \\
\frac{-1}{2\sqrt{n}}& \lambda &\frac{-1}{2} \\
0 &\frac{-1}{2} & \lambda \\
\end{array}\right)$.

\nt Then

$$ det(\lambda I - R(D_4^n))=\lambda (det(A))^n + 2n.\frac{1}{2\sqrt{n}}det
\left( \begin{array}{cccccc}
C& 0&0 &\ldots &0  \\
B&A&0 &\ldots&0   \\
B&0 &A&\ldots &0  \\
\vdots & \vdots &\vdots &\ddots &\vdots  \\
B&0& 0&\ldots &A   \\
\end{array} \right)_{3n\times 3n}. $$

\nt Now, by the straightforward computation we have the result.
	
\item[(ii)] It follows from Part (i).\quad\qed

\end{enumerate}

\section{Randi\'{c} energy of specific graphs with one edge deleted}

\nt In this section we obtain the randi\'{c} energy for certain graphs with one edge deleted.
 We need the following lemma:

\begin{lema}
Let $G=G_1\cup G_2\cup\ldots\cup G_m$. Then
$$RE(G)=Re(G_1)+ RE(G_2)+\ldots+ RE(G_m).$$
\end{lema}

\nt Here we state the  following easy  results:

\begin{lema}
\begin{enumerate}
\item[(i)] If $e\in E(P_n)$, then $RE(P_n-e)=RE(P_r)+RE(P_s)$, where $r+s=n$.

\item[(ii)] If $e\in E(C_n)$, ($n\geq 3$), then
$RE(C_n-e)=RE(P_n).$

\item[(iii)] Let $S_n$ be the star on $n$ vertices and $e\in E(S_n)$. Then for any $n\geq 3$,
$$RE(S_n-e)=RE(S_{n-1})=2.$$
\end{enumerate}
\end{lema}

\begin{teorem}
For $n\geq 2$,
\begin{itemize}
\item[(i)] Let $e$ be an edge of complete graph $K_n$. The Randi\'{c} characteristic polynomial of  $K_n-e$ is
$$RP(K_n-e,\lambda)=\lambda(\lambda -1)(\lambda +\frac{2}{n-1})(\lambda +\frac{1}{n-1})^{n-3}.$$
\item[(ii)] The Randi\'{c} energy of $K_n-e$ is
 $$RE(K_n-e)=2.$$
\end{itemize}
\end{teorem}
\noindent{\bf Proof.}
\begin{enumerate}
\item[(i)]
The Randi\'c matrix of $K_n-e$ is
$$\left( \begin{array}{cc}
0_{2\times 2}& \frac{1}{\sqrt{(n-1)(n-2)}}J_{2\times (n-2)} \\
\frac{1}{\sqrt{(n-1)(n-2)}}J_{(n-2)\times 2}&\frac{1}{n-1}(J-I)_{n-2}   \\
\end{array} \right).$$
Therefore
$$det(\lambda I -R(K_n-e))=det
\left( \begin{array}{cc}
\lambda I_2& \frac{-1}{\sqrt{(n-1)(n-2)}}J_{2\times (n-2)} \\
\frac{-1}{\sqrt{(n-1)(n-2)}}J_{(n-2)\times 2}&\frac{\lambda}{n-1}(J-I)_{n-2}   \\
\end{array} \right).
$$

Similar to the proof of Theorem \ref{bipartite}(i) we have the result.
\item[(ii)] It follows from Part (i).\quad\qed
\end{enumerate}

\begin{teorem}
For $m,n\neq 1$,
\begin{itemize}
\item[(i)] The Randi\'{c} characteristic polynomial of complete bipartite graph $K_{m,n}-e$ is
$$RP(K_{m,n}-e,\lambda)=\lambda^{m+n-4}(\lambda^2 -1)(\lambda^2 -\frac{1}{mn}).$$
\item[(ii)] The Randi\'{c} energy of $K_{m,n}-e$ is
 $$RE(K_{m,n}-e)=2+\frac{2}{\sqrt{mn}}.$$
\end{itemize}
\end{teorem}
\noindent{\bf Proof.}
\begin{enumerate}
\item[(i)]
The Randi\'c matrix of $K_{m,n}-e$ is
$$\left( \begin{array}{cc}
0_{m\times m}& A \\
A^t& 0_{n\times n}  \\
\end{array} \right),$$
where $A=\left( \begin{array}{cc}
\frac{1}{\sqrt{mn}}J_{(m-1)\times (n-1)}& \frac{1}{\sqrt{n(m-1)}}J_{(m-1)\times 1}\\
\frac{1}{\sqrt{m(n-1)}}J_{1\times (n-1)} & 0_{1\times 1}  \\
\end{array} \right)$.So
$$det(\lambda I -R(K_{m,n}-e))=\lambda I_m det (\lambda I_n - A^t\frac{1}{\lambda}I_m A)
.
$$

 Similar to the proof of Theorem \ref{bipartite}(i) we have the result.

\item[(ii)] It follows from Part (i).\quad\qed
\end{enumerate}



\begin{thebibliography}{99}

\bibitem{Bap} R. B. Bapat, S. Pati, Energy of a graph is never an odd integer, {\it Bull. Kerala Math.
Assoc.}, 1 (2004), 129--132.



\bibitem{Boz} \c S. B. Bozkurt, A. D. G\"ung\"or, I. Gutman, A. S. \c Cevik, Randi\'{c} matrix and Randi\'{c}
energy, {\it MATCH Commum. Math. Comput. Chem}. 64 (2010) 239--250.

\bibitem{Boz1} \c S. B. Bozkurt, A. D. G\"ung\"or, I. Gutman, Randi\'{c} spectral radius and Randi\'{c} energy,
{\it MATCH Commum. Math. Comput. Chem}. 64 (2010) 321--334.



\bibitem{Cve} D. Cvetkovi\'{c}, M. Doob, H. sachs, {\it Spectra of graphs} - Theory and Aplication, Academic Press, New York, 1980.

\bibitem{erdos}  P. Erd\"{o}s,
A. R\'{e}nyi, V.T.  S\'{o}s,  On a problem of graph theory, {\it  Studia Sci. Math. Hungar.}, 1,  215--235 (1966).




\bibitem{Gut} I. Gutman, The energy of a graph: Old and new results, in: A. Betten, A.Kohnert, R. Laue, A. Wassermannn (Eds.), {\it Algebraic Combinatorics and Applications}, Springer-Verlag, Berlin, 2001, pp. 196--211.

\bibitem{Gut1} I. Gutman, Topology and stability of conjugated hydrocarbons. The dependence of total $\pi$-electron energy on molecular topology, {\it J. Serb. Chem. Soc.} 70 (2005) 441--456.
    
    \bibitem{Gut3} I. Gutman, B. Furtula, \c S. B. Bozkurt, On Randi\'{c} energy, {\it Linear Algebra Appl.}, 442 (2014) 50--57.

\bibitem{Gut2} I. Gutman, X. Li, J. Zhang, Graph energy, in: M. Dehmer, F. Emmert-Streib
(Eds.), {\it Analysis of Complex Networks. From Biology to Linguistics}, Wiley-VCH,
Weinheim, 2009, pp. 145--174.

\bibitem {Maj} S. Majstorovi\'{c}, A. Klobu·car, I. Gutman, Selected topics from the theory of graph
energy: hypoenergetic graphs, in: D. Cvetkovi\'{c}, I. Gutman (Eds.), {\it Applications
of Graph Spectra}, Math. Inst., Belgrade, 2009, pp. 65--105.


\bibitem{Ran} M. Randi\'{c}, On characterization of molecular branching, {\it J. Amer. Chem. Soc}. 97 (1975) 6609--6615.


\bibitem{Roj} O. Rojo, L. Medina, Construction of bipartite graphs having the same Randi\'{c} energy, {\it MATCH Commun. Math. Comput. Chem.} 68 (2012) 805--814.





\end{thebibliography}
\end{document}